\documentclass[12pt]{article}
\usepackage{amsmath} 
\usepackage{amsfonts}
\usepackage{amssymb}
 \newcommand{\xn}{x}
  \newcommand{\ri}{r}
 \newcommand{\xii}{\xi}
  \newcommand{\si}{s}   
 \newcommand{\tii}{\zt}  
\usepackage{epsfig}
\usepackage{graphics}
 \newcommand{\Hfg}{H^{f,g}}
   \newcommand{\Blag}{Blagoveshchenski\v\i\ }
   \newcommand{\Phifg}{\Phi^{f,g}}
   \newcommand{\Psifg}{\Psi^{f,g}}
\newcommand{\intL}{\int_0^L}
\newcommand{\intt}{\int_0^t}
\newcommand{\intT}{\int_0^T}
\newcommand{\ints}{\int_0^s}

\newcommand{\zt}{\tau}
\newcommand{\zdia}{~~\rule{1mm}{2mm}\par\medskip}  

\newcommand{\zdiaform}{\mbox{~~\zdia}}  

\newcommand{\ZIN}{\infty}
\newcommand{\zProof}{{\bf\underbar{Proof}.}\ } 
  
\newcommand{\zg}{\gamma} 
\newcommand{\zaa}{\alpha} 
\newcommand{\zzr}{{\rm I\hskip-2.1pt R}}

\newcommand{\ZBI}{\bibitem}

\newcommand{\ZD}{\;\mbox{\rm d}}    
  
\newcommand{\ZEP}{\epsilon}

\newcommand{\ZLA}{\label} 
   
\newcommand{\zl}{\lambda}

\newcommand{\ZSI}{\sigma}

\newtheorem{Theorem}{Theorem}  
\newtheorem{Corollary}[Theorem]{Corollary}

\newtheorem{Remark}[Theorem]{Remark}
\newtheorem{Definition}[Theorem]{Definition} 
\newtheorem{Example}[Theorem]{Example}

\author{
L. Pandolfi\thanks{Dipartimento di Scienze Matematiche ``Giuseppe Luigi Lagrange'', Politecnico di Torino, Corso Duca degli Abruzzi 24, 10129 Torino, Italy (luciano.pandolfi@polito.it)}
}
\title{Identification of a space varying coefficient   of a linear viscoelastic string of Maxwell-Boltzman type\thanks{
This papers fits into the research program of the GNAMPA-INDAM and has been written in the framework of the   ``Groupement de Recherche en Contr\^ole des EDP entre la France et l'Italie (CONEDP-CNRS)''.}}

\begin{document}

\maketitle

\noindent {\bf\underline{Abstract}:} In this paper we solve the problem of the identification of a coefficient which appears in the model of  a distributed system with persistent memory encountered in linear viscoelasticity (and in diffusion processes with memory). The additional data used in the identification are subsumed in the input output map from the deformation to the traction on the boundary.  We extend a dynamical approach to identification introduced by Belishev in the case of purely elastic (memoryless) bodies and based on a special equation due to  \Blag\!. So, in particular, we extend \Blag equation   to our class of systems with persistent memory.

\medskip


\smallskip 

\section{Introduction}

A noticeable property of viscoelastic materials is that the reaction to applied stimuli persists in time and when the memory is long a commonly used linear model for a viscoelastic string    is (see~\cite{Kolski,Pipkin}):
\begin{equation}
\ZLA{eq:modello}
w_t=\intt N(t-s)\mathcal{L} w(s)\ZD s\,,\qquad \mathcal{L}=\Delta+q(x)\,,\quad x\in (0,L)
\end{equation}
where $w=w(x,t)$ represents the displacement at time $t$ and position $x$ and $N(t)$ is the \emph{relaxation kernel.} 
Note that if $N(t)\equiv 1$   then Eq.~(\ref{eq:modello}) is the integrated version of the string equation
\begin{equation}
\ZLA{eq:wave}
w_{tt}=\mathcal{L} w \,,\qquad x\in (0,L)\,.
\end{equation}
 Equation~(\ref{eq:modello}) 	is also encountered in diffusion processes (when the material has a complex molecular structure)  and thermodynamics  and in this contest it was first introduced by Maxwell~(\cite{Maxwell}, see also~\cite{Cattaneo}) in the special case $N(t)=e^{-a t}$, $a>0$.

An important problem widely studied for the wave equation~(\ref{eq:wave}) is 
the identification of the unknown coefficient $q(x)$
  on the basis of experimental measurements of the deformation and traction on a part of the boundary:         a deformation $f(t)$ is applied to the boundary, for example at $x=0$, and the resulting traction $\ZSI^f(t)$ at the same part of the boundary is measured. In this way we construct the transformation $f \mapsto\ZSI^f$, defined in suitable spaces, and the coefficient $q(x)$ 
  in Eq.~(\ref{eq:wave}) can be reconstructed once this transformation is known.

We sum up: we assume that Eq.~(\ref{eq:modello})  is supplemented with the following initial/boundary conditions:
\begin{equation}
\ZLA{eq:condiINIboundary}
w(x,0)=0\,, \quad w(0,t)=f(t)\,,\qquad w(L,t)=0 \,.
\end{equation}
The corresponding solution is denoted $w^f=w^f(x,t)$. The measurement is
\begin{equation}
\ZLA{eq:BoundaryMEASU}
y^f(t)= w^f_x(0,t)
\end{equation}
and the goal is the reconstruction of the coefficient $q(x)$ from the map $f\mapsto y^f$.

 \begin{Remark}{\rm
\begin{itemize}
\item We must keep in mind that:
\begin{equation}\ZLA{eq:CondiDERIV}
\left\{\begin{array}{l}
 \mbox{the condition}\quad   
w_t(x,0)=0\quad  \mbox{follows from Eq.~(\ref{eq:modello})} \\
\mbox{while $w_t(x,0)=0$  has to be added when working with Eq.~(\ref{eq:wave})}\\
\end{array}\right.
\end{equation}

\item When $f$ is merely of class $L^2$, we cannot hope for a smooth solution $w$ of Eq.~(\ref{eq:modello}). The definition of \emph{mild solution} is in Sect.~\ref{sectOBERCmodel}.
\item The boundary observation $y^f(t)$ is square integrable   if the boundary deformation is smooth, 
see Theorem~\ref{TeoSEC2RegoLAR} below. 
\item The traction exerted on the string by its support   is 
\[\ZSI^f(t)=-\intt N(t-s)y^f(s)\ZD s\,.
\] In practice, this is the function that it is measured. 
Thanks to the  regularity of the relaxation kernel $N(t)$ stated below, we can compute $y^f(t)$ from the Volterra integral equation
\[
y^f(t)+\intt N'(t-s)y^f(s)\ZD s=-\frac{\ZD}{\ZD t}\ZSI^f(t) \,.\zdiaform
\]

\end{itemize}
}\end{Remark}

  Few comments on the literature: among the several papers on this dynamical method for the reconstruction of $q(x)$ in a wave equation we mention the very readable papers~\cite{AvdoninBeliIVANOV,BELIwaves}. 
  See~\cite{BeliSurvey1,BeliSurvey2} for further advances and extensions.    
 We note that the dynamical approach to the identification of $q(x)$ has been used also in the books~\cite{Romanov,Romanov1}. See for example~\cite{Morassi} for an application of the methods in~\cite{Romanov,Romanov1} to a concrete problem 
 and~\cite{RomanovMEMORY} and references therein for several inverse problems and methods for viscoelastic materials. 
Finally, we mention a related approach, which assume that we know the \emph{spectral data,}    i.e.   the pairs $\{\zl_n,\phi_n'(0)\}$ where   $\{\phi_n(x)\}$ is a complete orthonormalized sequence of eigenfunctions and $\zl_n$ is the eigenvalue of $\phi_n$. See~\cite{Katchalov} for this approach and note that the spectral data can be derived from the map $f\mapsto y^f$, see for example~\cite{AvdLenProto}.

The organization of the paper is as follows: the crucial ingredient in the algorithm used to identify $q(x)$ for the \emph{wave} equation is a special equation introduced by 
\Blag (see~\cite{BLAGO}. In this paper the equation is written in the case a boundary input is a delta function). This equation is extended to system~(\ref{eq:modello}) in sect.~\ref{sec:Blag} and we believe that this is the main contribution of the paper. Based on this new version of \Blag equation, the algorithm for the identification of $q(x)$ for the string equation can be extended to  system~(\ref{eq:modello}). This is done in
 Sect.~\ref{sect:identifQ}.
 
\subsection{\ZLA{subsec:ASSUnota}Assumptions and notation}
The standing assumptions in this paper are:
\begin{itemize}

\item We assume that $N(t)$ is of class $C^3$, with $N(0)>0$, which accounts for the iperbolic nature of the process~(\ref{eq:modello}). For simplicity of notations, we introduce the nonrestrictive assumption that the time scale has been normalized so that
\[
N(0)=1\,.
\]
 \item the function $q(x)$ is continuous.
\item
The (unbounded) map $f\mapsto y^f=w^f_x(0,\cdot)$: $L^2(0,T)\mapsto L^2(0,T)$ is denoted $R_T$:
\[
R_T f=y^f
\]
and it is called the \emph{response operator.  We assume that the map $R_T$ is known.}

\end{itemize}

We list few properties which follows from the assumption on the relaxation kernel $N(t)$:
\begin{enumerate}
\item\ZLA{itemCONTIN} Let $T>0$. for every $f\in L^2(0,T)$ the (mild) solution $w^f$ of Eq.~(\ref{eq:modello}) takes values in $C([0,T];L^2(0,L))$ 
and $f\mapsto w^f$ is a linear and continuous transformation in these spaces (see Sect.~\ref{sectOBERCmodel} for the definition of the mild solutions);
\item the velocity of wave propagation in the viscoelastic string discribed by~(\ref{eq:modello}) is $N(0)=1$, as in the purely elastic string described by the equation~(\ref{eq:wave});
\item \ZLA{itemCONTROLL} system~(\ref{eq:modello}) is controllable in time $L$, i.e. \emph{for every $\xi(x)\in L^2(0,L)$ there exists $f\in L^2(0,T)$, $T=L$ such that $w^f(x,T)=\xi(x)$.} Here $w^f(x,t)$ is the mild solution of Eq.~(\ref{eq:modello}) and the actual meaning of controllability (which can be equivalently stated in different ways) is as follows: the transformation $f(\cdot)\mapsto w^f(\cdot,T)$ is surjective from $L^2(0,L)$ to itself. In terms of the sine Fourier expansion, controllability can be stated as follows. We expand  $\xi(x)=\sum_{n=1}^{+\ZIN} \xi_n\sin(\pi/L) nt$, $w(x,t)=\sum_{n=1}^{+\ZIN} w^f_n(t)\sin(\pi/L) nt$. Controllability is the property that there  exists $f\in L^2(0,T)$ such that $w_n^f(T)=\xi_n$ for every index $n$.
\end{enumerate}
See~\cite[Theorem~3.4]{PandIEOT} for the case $q=0$ and~\cite[Sect.~6.1]{PandCINA} for $q\neq 0$. See~\cite{Belleni,DEschGrimmer,PandLIBRO} for the propagation speed ($q(x)=0$ in these references,  but the proof is easily adapted).

\paragraph{Notations}
\begin{itemize}
  \item We shall have to consider two intervals of the same length $T$ (due to the fact that the velocity of wave propagation is $1$). The first interval is contained  in the time axis and the second one in the space axis. When needed for clarity, we keep them distinct notationally as follows: on the time axis the interval is denoted $ [0, T]$ while in the  space axis it is denoted $ [0,  T_X ]$ but $  T_X=T$.  
  
  With this notations, the controllability property stated in item~\ref{itemCONTROLL} can be expressed more clearly:
  the transformation $f(\cdot)\mapsto w^f(\cdot,T)$ is surjective from $L^2(0,  L)$ to $L^2(0, L_X)$, so stressing the fact that these spaces have different roles, even if $  L= L_X $.  
\item Let $T>0$ be fixed. We denote $\Lambda_T$ the map
\[
\Lambda_Tf=w^f(\cdot,T)\,.
\]
The properties   stated in item~\ref{itemCONTIN}  above imply   $\Lambda_T \in \mathcal{L}(L^2(0,T);L^2(0,L))$ for every $T>0$ and controllability implies that $\Lambda_T$ is surjective for $T$  ``large''  $T\geq L$ (both in the case of the wave equation and of Eq.~(\ref{eq:modello}), see~\cite{PandCINA,PandLIBRO}).

\item The map $C_T=\Lambda_T^*\Lambda_T\in \mathcal{L}\left (L^2(0,T)\right )$   is the \emph{controllability operator} for the controllability of the deformation $w$. In the contest of inverse problems it is called the \emph{connecting operator.}

The operator $C_T$ is selfadjoint positive and boundedly invertible in $L^2(0,T) $  if $T\geq L$, thanks to controllability.
 
 \item
A prime is used to denote the derivative of a function of only one variable. In the case of functions of more variables, the derivative is denoted either with the corresponding variable put at the index, or with a comma followed by the number which identifies the position of the variable, as common in engineering literature. So, for example
\[
w_{,1}(x,t)=w_x(x,t)\,,\qquad w_{,11}(x,t)=w_{xx}(x,t)\,,\qquad w_{,2}(x,t)=w_{t}(x,t)\, \dots
\] 

 \item The function $N'(t)$ has a special role (the wave equation corresponds to $N'(t)=0$). We put
 \[
N_1(t)=N'(t)\,. 
 \]

 We shall encounter several Volterra integral equations whose kernel is $N_1(t)$, i.e. equations of the form
 \begin{equation}\ZLA{eq:VolteGEN}
v(t)+\intt N_1(t-s) v(s)\ZD s=F(t) \,.
 \end{equation}
 The solution $v(t)$ of this equation is
 \begin{equation}\ZLA{eq:VolteGENsolu}
v(t)=F(t)-\intt R(t-s) F(s)\ZD s 
 \end{equation}
 where $R(t)$ is the \emph{resolvent kernel} of $N_1(t)$, i.e. the unique solution of the equation
\begin{equation}\ZLA{eq:Risolvente}
R(t)+\intt N_1(t-s)R(s)\ZD s=N_1(t)\,.
\end{equation}
We assumed that $N_1(t)\in C^3$ so that we have also $R(t)\in C^3$.
\end{itemize}

\section{\ZLA{sectOBERCmodel}Preliminaries}

We shall repeatedly use    ``MacCamy trick'', as in~\cite[Sect.~2.2]{PandLIBRO}. MacCamy trick is a formal manipulation, which goes as follows: first we compute the derivatives of both the sides in~(\ref{eq:modello}) and we get
 \[
w''(x,t)=(\Delta+q)w(x,t)+\intt N'(t-s) (\Delta+q)w(x,s)\ZD s\,.
 \]
  For every fixed value of $x$,
 this is a Volterra integral equation of the form~(\ref{eq:VolteGEN}) in the unknown $(\Delta+q)w(x,t)$. We  solve  it  using formula~(\ref{eq:VolteGENsolu}) and, using $R(t)\in C^3$, we integrate $\intt R(t-s)w''(x,s)\ZD s$ by parts. 
 Using $w(x,0)=0$, $w_t(x,0)=0$ we get
 \begin{align*}
w''(x,t)&=R(0) w'(x,t)+R(0) w(x,t)+\intt R''(t-s) w(x,s)\ZD s + \\
&+ \left (w_{xx}(x,t) +q(x)w(x,t)\right )\,.
 \end{align*}
 We introduce the notations
 \begin{equation}\ZLA{eq:NOtazioni}
 \zg=\frac{R(0)}{2}\,,\qquad \zaa=R'(0)+R^2(0)/4\,, \qquad 
K(t)=e^{-\zg t}R''(t) 
 \end{equation}
and we replace $w(x,t)=e^{\zg t}W(x,t)$.
It turns out that the function $W(x,t)$ solves
\begin{equation}\ZLA{eq:DOPOmaCCamy}\begin{array}{l}
\displaystyle 
W''(x,t)= \Delta W(x,t)+\left (q(x)+\zaa\right ) W(x,t)+ \intt K(t-s) W(x,s)\ZD s\,,\\[3mm]
W(x,0)=0\,,\quad W_t(x,0)=0\,,\qquad W(0,t)=e^{-\zg t} f(t)\,,\quad W(L,t)=0\,.
\end{array}   
\end{equation}


Eq.~(\ref{eq:DOPOmaCCamy}) is a (memoryless) string equation ``perturbed'' by the term
\begin{equation}\ZLA{eq:LaFORMperty}
F(x,t)=\left (q(x)+\zaa\right )W(x,t)+\intt K(t-s) W(x,s)\ZD s\,.
\end{equation}
 By definition, a \emph{mild solution} of Eq.~(\ref{eq:modello}) is a function $w(x,t)$ such that $W(x,t)=e^{-\zg t}w(x,t)$ is a mild solution of the perturbed wave equation~(\ref{eq:DOPOmaCCamy}). So,
the functions $w(x,t)$ and $W(x,t)$ have the same regularity properties.

We shall use the following known fact (see~\cite{TikoSAMARSKI}): the solutions of the wave equation
\[
W_{tt}(x,t)=W_{xx}(x,t)+F(x,t)\,,\quad x>0\,,\qquad 
  \left\{
 \begin{array}
 {l}
 W(\xn,0)=0\,, \\
 W(0,t)=  g(t) 
 \end{array}
 \right.
\]
are given by the following formula:
 
\begin{equation}\ZLA{FormDALA}
\begin{array}{l}
\displaystyle  W(x,t)= g(t-x)+\frac{1}{2}\int _{D(x,t)} F(\xii,\tii)\ZD\xii\,\ZD\tii=\\[3mm]
\displaystyle  = g(t-x)+\frac{1}{2}\intt \int _{|x-t+\tii|}^{x+t-\tii} F(\xii,\tii)\ZD\xii\,\ZD\tii\,.
\end{array}
\end{equation}
This same formula represents the solutions also if $x\in (0,L)$   with boundary condition $W(L,t)=0$, provided that we confine ourselves to consider solely $t\in (0,T)$,  $T\leq L$ i.e. before a reflection from the right part of the boundary appears. 
So, when $T\leq L$ we can apply formula~(\ref{FormDALA}) to our problem~(\ref{eq:modello}) with the initial/boundary conditions~(\ref{eq:condiINIboundary}).
In this case, $g(t)=e^{-\zg t} f(t)$ and $F$ is the function in~(\ref{eq:LaFORMperty}) so that we get   a Volterra integral equation in $L^2(0,L)$ for the unknown function $W(x,t)$, on the interval of time $[0,T]\subseteq [0,L]$ and   we can state that  $w\in C([0,T];L^2(0,L))$, $T\leq L$, is a \emph{mild solution} of Eq.~(\ref{eq:modello}) with conditions~(\ref{eq:condiINIboundary})when $W$ is a continuous solution of this Volterra integral equation.

We need   information on  the solutions of Eq.~(\ref{eq:modello}) with conditions~(\ref{eq:condiINIboundary}), which are contained, more or less implicitly, in several papers, see for example~\cite{Belleni,DEschGrimmer,PandLIBRO}. 

In order not to interrupt the presentations, the proofs of these instrumental results  are in the appendix.

The first information we need is on the regularity of the solutions of Eq.~(\ref{eq:modello}) in the style of~\cite[Theorem~2.4]{PandLIBRO}. We need this information only for the case of the viscoelastic string, and only when $T\leq L$. So we present an elementary derivation of the following result:
\begin{Theorem}\ZLA{Teo:RegoSOLUZ}
Let $T\leq L$ and let $f\in\mathcal{D}(0,T)$. Then $w_{tt}(x,t)$ and $w_{xx}(x,t)$ are continuous functions and equality~(\ref{eq:modello}) holds on $(0,L)\times(0,T)$.
\end{Theorem} 

We noted already that disturbances propagates with the same speed as in the case of the corresponding wave equation.  A consequence is as follows: we fix a time instant $T\in(0,L)$. We consider control 
system~(\ref{eq:modello}),~(\ref{eq:condiINIboundary}) but with $x\in(0,T_X)$ and   the condition $w(T_X,t)=0$ (the notations $T$ and $T_X $ are as in 
the subsection~\ref{subsec:ASSUnota} ``Assumptions and notation''. So, $ T_X=T$).  The condition $w(x,t)=0$ for $x>T$
  is automatically satisfied since the propagation speed is precisely $1$. So, we have a system on an interval $ (0,T_X)$ which   is controllable at time $ T=T_X$ (see~\cite{PandCINA}) and:
  \begin{Corollary}\ZLA{CoroREGOL}
 Let $T\leq L$.  The following properties hold:
\begin{itemize}  
  \item
  for every $\xi\in L^2(0,T_X)$ (we recall: $T_X=T$) there exists $f\in L^2(0,T )$ such that $w^f(x,T)=\xi(x)$. 
  \item
  If we confine ourselves to use controls $f\in \mathcal{D}(0,T)$   then:
  \begin{itemize}
  \item the reachable set $\{w^f(x,T)\,\ f\in \mathcal{D}(0,T)\}$ is dense in $L^2(0,T_X)$;
 \item the support of $f\in\mathcal{D}(0,T)$ is properly contained in $[0,T -\ZEP]$ for some $\ZEP>0$ so that for $x\in (T_X,L]$ we have $w^f(x,T)=w^f_x(x,T)=0$.
   \end{itemize} 
  \end{itemize}   
\end{Corollary}

Now we state the following result (see  the Appendix for the proof):
 \begin{Theorem}\ZLA{TeoSEC2RegoLAR}Let $T\leq L$. The following properties hold:
 \begin{enumerate}
\item\ZLA{TeoSEC2RegoLARitem1} Let   the target $\xi(x)$ be continuous, $\xi\in C([0,T_X])$ (we recall the notations: $T_X=T$). Then, there exists a unique \emph{steering control,}  i.e. a control $f$ such that   $w^f(x,T)=\xi(x)$,    and this control is continuous.

\item\ZLA{TeoSEC2RegoLARitem1BIS}
If the target $\xi(x)$ is continuous   then
  \begin{equation}\ZLA{eq:delLIMITaZERO}
 \lim _{x\to T ^-} w(x,T )=\lim _{t\to 0^+} f(t)\,,\quad \mbox{and, when $T<L$,}\quad  \lim _{x\to T ^+} w(x,T )=0\,.  
  \end{equation}
\item\ZLA{TeoSEC2RegoLARitem2}
The response operator $(R_T f)(t)=y^f(t)$ is unbounded in $L^2(0,T)$. Its domain is $H^1(0,T)$.
\item\ZLA{TeoSEC2RegoLARitem4} If $w^f(x,T)=\xi(x)\in H^1(0,T_X) $ then $f\in H^1(0,T)$.
\end{enumerate}
\end{Theorem}

 See~\cite{pandTRIULZI} for a   weaker form of the last statement. 

\section{\ZLA{sec:Blag}\Blag equation for system~(\ref{eq:modello})}

For clarity we outline the derivation of \Blag equation for the string equation~(\ref{eq:wave}) with the initial and boundary conditions in~(\ref{eq:condiINIboundary})-(\ref{eq:CondiDERIV}). 

We fix any $T>0$ and we consider two smoothly varying boundary deformations  $f$ and $g$ ($f,\, g\in H^2(0,T)$ is enough to justify the following computations).  
Let $\Hfg(s,t)$ be the function defined on $[0,T]\times[0,T]$ by
\[
\Hfg(s,t)=\int_0^L w^f(x,t)w^g(x,s)\ZD x\,.
\]
The initial conditions of Eq.~(\ref{eq:wave}) are zero so that
\begin{equation}\ZLA{eq:COndiINIBOUnddiH}
\Hfg(0,t)=\Hfg(s,0)=0\,,\qquad \Hfg_t(s,0)=\Hfg_s(0,t)=0\,.
\end{equation}
An integration by parts gives
\begin{align*}
&\Hfg_{tt}(s,t)=\int_0^L w^f_{tt}(x,t)w^g(x,s)\ZD x=\\
&=\int_0^L w^f_{xx}(x,t)w^g(x,s)\ZD x+\int_0^L w^f(x,t)q(x)w^g(x,s)\ZD x=\\
&=w_x^f(L,t)w^g(L,s)-w_x^f(0,t)w^g(0,s)-w^f(L,t)w_x^g(L,s)+w^f(0,t)w_x^g(0,s)+\\
&+
\int_0^L w^f (x,t)\left [ w^g_{xx}(x,s) +  q(x)w^g(x,s)\right ]\ZD x 
 =\Hfg_{ss}(s,t)+\Phifg(s,t)\,,\\
& \Phifg(s,t)=f(t)y^g(s)-y^f(t)g(s)=
f(t)(R_Tg)(s)-(R_Tf)(t)g(s)
\,.
\end{align*}
Our assumption is that the response operator $R_T$ is known so that 
the function $\Phifg(s,t)$ is known and the function $\Hfg(s,t)$ can be computed (for every $f$, $g$ smooth) from the equation $\Hfg_{tt}(s,t)=\Hfg_{ss}(s,t)+\Phifg(s,t)$ with the boundary 
conditions~(\ref{eq:COndiINIBOUnddiH}) (one of the boundary conditions on the derivative is redundant of course, but we should not worry about that, since the existence of $\Hfg(s,t)$ is clear from its definition). 

In particular, for every $t$ and every $f$ and $g$ smooth, we can compute
\[
\Hfg(t,t)=\int_0^L \left ( \Lambda_tf\right )(x)\left ( \Lambda_t g\right )(x)\ZD x=\langle \Lambda_t^*\Lambda_t f,g\rangle_{L^2(0,T)}.
\]
But, the operator $f\mapsto \Lambda_t^*\Lambda_t f$ is continuous on $L^2(0,t)$ for every $t$ and so \Blag equation shows that the controllability operator $f\mapsto \Lambda_t^*\Lambda_t f$ can be computed from the response map even if this map is defined only for smooth functions.  This is the key for the reconstruction of the coefficient $q(x)$ in the string equation.

In this section we prove that an analogous of the \Blag equation can be derived for the system described by~(\ref{eq:modello}).

The definition of $\Hfg(t,s)$ is the same as above,
\begin{equation}\ZLA{eq:defiHfg}
\Hfg(s,t)=\int_0^L w^f(x,t)w^g(x,s)\ZD x
\end{equation}
but now $w^f$ and $w^g$ solve~(\ref{eq:modello}) with the conditions~(\ref{eq:condiINIboundary}). Using the fact that $w_t(x,0)=0$  (see the observation 
at the first line of~(\ref{eq:CondiDERIV})) we see that  also in the case of Eq.~(\ref{eq:modello}) the function $\Hfg(s,t)$ satisfies the  
conditions~(\ref{eq:COndiINIBOUnddiH}). Furthermore,
 
\begin{align*}
&\nonumber\frac{\partial}{\partial t}\Hfg (s,t)=\intL w^f_t(x,t)w^g(x,s)\ZD x= \\
&\intL\left [\intt N(t-r)\left [w^f_{xx}(x,r)+q(x)w^f(x,r)\right ]\ZD r\right ]w^g(x,s)\ZD x=
\\
&\nonumber=
\intt N(t-r)\intL w^f(x,r)q(x)w^g(x,s)\ZD x\,\ZD r+\\
&+\intt N(t-r)\intL w^f(x,r)w^g_{xx}(x,s)\ZD x\,\ZD r+\Phifg(s,t)=\\
& =\Phifg(s,t)+\intt N(t-\xi)\intL w^f(x,\xi)\left [q(x)w^g(x,s)+w_{xx}^g(x,s)\right ]\ZD x\,\ZD\xi
\end{align*}
where now
\[
\Phifg(s,t)=\intt N(t-r)\left [ f(r)y^g(s)-y^f(r)g(s)\right ]\ZD r\,.
\]
We compute
\begin{align*}
&\int_0^s N(s-\ZSI)  \frac{\partial}{\partial t}\Hfg(\ZSI,t)\ZD\ZSI=\int_0^s N(s-\ZSI)\Phifg(\ZSI,t)\ZD\ZSI+\\
&+\int_0^s N(s-\ZSI)\intt N(t-\xi)\intL w^f(x,\xi)\left [ q(x) w^g(x,\ZSI)+w^g_{xx}(x,\ZSI)\right ]\ZD x\,\ZD \xi\,\ZD\ZSI=\\
&=\intt N(t-\xi)\intL w^f(x,\xi)\int_0^s N(s-\ZSI)\left [q(x)w^g(x,\ZSI)+w_{xx}^g(x,\ZSI)\right ]\ZD\ZSI\,  \ZD x\, \ZD\xi+\\&+\int_0^s N(s-\ZSI)\Phifg(\ZSI,t)\ZD\ZSI=\\
&=\intt N(t-\xi )\frac{\partial}{\partial s}\Hfg (s,\xi)\ZD\xi+\int_0^s N(s-\ZSI)\Phifg(\ZSI,t)\ZD\ZSI \,.
\end{align*}
Hence we have
\begin{equation}
\ZLA{Eq:BlagoPRELI}
\left\{\begin{array}{l}\displaystyle
\int_0^s N(s-\ZSI)\frac{\partial}{\partial t} \Hfg (\ZSI,t)\ZD \ZSI=\intt N(t-\xi)\frac{\partial}{\partial s}\Hfg(s,\xi)\ZD \xi+\Psifg(s,t )\,,\\
\displaystyle \Psifg(s,t )=\int_0^s N(s-\ZSI)\Phifg( \ZSI,t)\ZD\ZSI\,.
\end{array}\right.
\end{equation}
The function $\Psifg( \si,t)$ is known, computed from the response operator so that $\Hfg(\si,t)$ solves
\[
\frac{\partial}{\partial  t}\ints N(s-\ZSI)\Hfg(\ZSI,t)\ZD\ZSI=\frac{\partial}{\partial s} \intt N(t-\xi)\Hfg(s,\xi)\ZD\xi+\Psifg(s,t)\,.
\]
 
We compute the derivative of both the sides respect to $s$ and we get (we recall the notation $N_1(t)=N'(t)$):
\begin{align}
\nonumber&=
\left [\frac{\partial}{\partial t}\Hfg(s,t)\right ]+\ints N_1(s-\ZSI) \left [ \frac{\partial}{\partial t}\Hfg(\ZSI,t)  \right ]\ZD\ZSI=
\\
\ZLA{eq:nuovoPAG4}
&=\frac{\partial^2}{\partial s^2}\intt N(t-\xi)\Hfg(s,\xi)\ZD\xi+ 
 \frac{\partial}{\partial s}\Psifg(s,t)\,.
\end{align}
For every fixed value of the parameter $t$, we have a Volterra integral equation in the variable $s$. 
Then, using~(\ref{eq:VolteGENsolu}) as in Sect.~\ref{sectOBERCmodel},  we have
\begin{align}\nonumber
&\frac{\partial}{\partial t}\Hfg(s,t)=\frac{\partial^2}{\partial s^2}\intt N(t-\xi)\Hfg(s,\xi)\ZD\xi -\\
\ZLA{eq:BlagoMEM}&
 -\ints R(s-\nu)\left [
 \intt N(t-\xi )\frac{\partial^2}{\partial \nu^2}\Hfg(\nu,\xi)\ZD\xi\,\ZD\nu
\right ]+F(s,t)
\end{align}
with  
\[
 F(s,t)=\frac{\partial}{\partial s}\Psifg(s,t)-\ints R(s-\nu)\frac{\partial}{\partial \nu}\Psifg( \nu,t )\ZD\nu\,.
\]
 
Equality~(\ref{eq:BlagoMEM})   holds for every $t>0$ and $s>0$ and  it
  is  the \emph{\Blag  equation} for the system~(\ref{eq:modello}).  In this equation, the affine term $F(s,t)=F^{f,g}(s,t)$ can be computed from $f$ and $g$ via the response operator. Furthermore we know that $\Hfg(s,t)$ satisfies the conditions~(\ref{eq:COndiINIBOUnddiH}).

  When the inputs $f$ and $g$ are smooth,    it is possible to rewrite Eq.~(\ref{eq:BlagoMEM})
  with conditions~(\ref{eq:COndiINIBOUnddiH})
  in the form of a Volterra integral equation (see Eq.~(\ref{eq:intePERblago}) below). By definition, the \emph{mild solutions} of~(\ref{eq:BlagoMEM}) are the solution of the Volterra integral equation~(\ref{eq:intePERblago}).

\begin{Remark}
{\rm
If $f$, $g$ are of class $H^1$ then
\begin{align*}
\Psifg(s,t)&=\left [
\ints N(s-\ZSI)y^g(\ZSI)\ZD\ZSI
\right ]\intt N(t-r) f(r)\ZD r-\\
&-\left [
\ints N(s-\ZSI)g(\ZSI)\ZD\ZSI
\right ]\intt N(t-r)y^f(r)\ZD r
\end{align*}
is continuous in $(t,s)$ and of class $H^1_{\rm loc}$.\zdia
}
\end{Remark}
In the course of the proof of the next result,  we give a definition of the mild solutions of \Blag equation and we prove:  
\begin{Theorem}
Let $T>0$ be fixed and let $f$, $g\in H^1(0,T)$. Then,
Eq.~(\ref{eq:BlagoMEM}) with data~(\ref{eq:COndiINIBOUnddiH}) admits a unique \emph{mild solution} which is continuous on $[0,T]\times [0,T]$.
\end{Theorem}
\zProof 
We perform formal computations  which can be justified if the data are sufficiently smooth. \emph{We derive  an integral equation whose solutions are by definition the mild solutions of the  \Blag  equation~(\ref{eq:BlagoMEM}).} Then we prove existence and unicity of solutions of this integral equation.

We must use several changes of variables. So, in this computation it is convenient to  use the comma notation for the derivative: 
$H_{,1}$ and $H_{,2}$ denote respectively the derivative respect to the first and the second variable. The apex denotes the derivative of a function of one variable.

We integrate by parts  the last integral in~(\ref{eq:BlagoMEM}) and we use $H(0,\xi)=H_1(0,\xi)=0$. We get

\begin{align*}
  \Hfg_{,2}(s,t)&= \intt N(t-\ri)\Hfg_{,11}(s,\ri)\ZD\ri -
 \intt N(t-\ri)\left [ R(0)\Hfg_{,1}(s,\ri)+R'(0)\Hfg(s,\ri)\right] \ZD r-\\
&-\intt N(t-\ri)\int_0^{s}R''(s-\nu)\Hfg(\nu,\ri)\ZD\nu\,\ZD\ri+F(s,t)\,.
\end{align*} 
We replace $\Hfg(s,t)$ with
\[
\Hfg(s,t)=e^{\zg s}V(s,t)\,,\qquad \zg=R(0)/2\,.
\]
This removes the derivative respect to the first variable in the second integral of the first line and gives the following equation for $V$:
 
\begin{align*}
&V_{,2}(s,t)= \intt N(t-\ri)\left [\zaa V(s,\ri)+V_{,11}(s,\ri)\right ]\ZD\ri-\\
&-\intt N(t-\ri)\int_0^{s}R_2(s-\nu)V(\nu,\ri)\ZD\nu\,\ZD \ri+F(s,t) e^{-\zg s}
\end{align*} 
where
\begin{equation}\ZLA{eq:LeDEFalphaECC}
\zaa=-R'(0)-R^2(0)/4\,,\qquad R_2(\xii)=R''(\xii)e^{-\zg \xii}=R''(\xii)e^{-\left (R(0)/2\right )\xii}\,:
\end{equation}
We compute the derivative of both the sides respect to the second variable $t$ and we get
\begin{align*}
 &\left [V_{,11}(s,t)+\zaa V(s,t)-\int_0^{s} R_2(s-\nu)V(\nu,t)\ZD\nu\right ]+\\
 &+\intt N_1(t-\ri) 
 \left[  V_{,11}(s,\ri)+\zaa V(s,\ri)-\int_0^{s} R_2(s-\nu)V(\nu,\ri)\ZD\nu
 \right ]\ZD\ri=\\
&=V_{,22}(s,t) +F_{2}(s,t)\,,\qquad F_2(s,t)=-\frac{\ZD}{\ZD t}e^{-\zg t}F_{,2}(s,t)\,.
\end{align*} 
We use again that $R(t)$ is the resolvent kernel of $N_1(t)$. We solve for the unknown  bracket in the left side, this time respect to the variable $t$, for any fixed $s$. We get:
\begin{align*}
&\left [V_{,11}(s,t)+\zaa V(s,t)-\int_0^{s} R_2(s-\nu)V(\nu,t)\ZD\nu\right ]=V_{,22}(s,t)-\intt R(t-\ZSI)V_{,22}(s,\ZSI)\ZD\ZSI+\\
&+ F_{2}(s,t)-\intt R(t-\ZSI)F_{2}(s,\ZSI)\ZD\ZSI 
\\
&=
V_{,22}(s,t)-\left [
R(0)V_{,2}(s,t)+R'(0)V(s,t)+\intt R''(t-\ZSI)V(s,\ZSI)\ZD\ZSI 
\right ]+\hat G(s,t)\,,\\
&\hat G(s,t)=F_{2}(s,t)-\intt R(t-\ZSI)F_{2}(s,\ZSI)\ZD\ZSI\,.
\end{align*}
 
Finally, we replace 
\[
V(s,t)=e^{\zg t}W(s,t)\,,\qquad \zg=   R(0)/2 
\]
and we remove the derivative $V_{,2}(s,t)$ from the right hand side. We take into account the definition of $\zaa$, $\zaa=\zaa=-R'(0)-R^2(0)/4$ in~(\ref{eq:LeDEFalphaECC}) and   we get  

\begin{align*}
&W_{,22}(s,t)=W_{,11}(s,t)+\intt R_2(t-\zt)W(s,\zt)\ZD\zt-\\
&-\int_0^{s} R_2(s-\nu)W(\nu,t)\ZD\nu+G(s,t)\,,\qquad 
 G(s,t)=-e^{-\left(R(0)/2\right) t}\hat G(s,t) 
\end{align*}
 (the terms with the factor $\zaa$ cancel out).  
 This is an equation on the quarter plane $t>0$, $s>0$ with ``initial'' and ``boundary'' conditions $W(s,0)=W_t(s,0)=0$,  $W(0,t)=0$.  
We apply   formula~(\ref{FormDALA}) (now with the boundary input put equal zero) and we get the following integral equation for $W(\xn,t)$:
\begin{align}
\nonumber
W(s,t)&= \frac{1}{2} \int _{D(s,t)} \left [ 
\int_0^{\tii}R_2(\tii-\nu)W(\xi,\nu)\ZD\nu-\right.\\
\ZLA{eq:intePERblago}&\left.-\int_0^{\xii}R_2(\xii-\nu)W(\nu,\tii)\ZD\nu
\right ]\ZD\xii\, \ZD\tii+
\frac{1}{2}
\int _{D(s,t)} 
G(\xii,\tii)
 \ZD\xii\,\ZD\tii\,.
\end{align}
This is the integral version of the \Blag equation and 
\emph{by definition the function $\Hfg(\si,t)$ is the mild solution of~(\ref{eq:BlagoMEM}) with data~(\ref{eq:COndiINIBOUnddiH}) when $W(s,t)$ solves the integral equation~(\ref{eq:intePERblago}).}

Note that if $f$, $g$ are of class $H^1$ then the last integral is a continuous function of $(s,t)$.

Now we prove that Eq.~(\ref{eq:intePERblago}) admits a unique continuous solution. 
It is convenient to make a last transformation: we introduce
\[
 W(s,t)=  e^{\ZSI(s+t)}Y(s,t)
\]
so that $Y(s,t)$ solves  
\begin{align}
\nonumber Y(s,t)
 &
 = \frac{1}{2} \int _{D(s,t)} \left [ 
  \int_0^{\tii}\left ( R_2(\tii-\nu)e^{-\ZSI\left ((s-\xi)+(t-\nu)\right )}\right)Y(\xii,\nu)\ZD\nu -\right.
 \\
\nonumber &
   \left.-
 \int_0^{\xii}\left (R_2(\xii-\nu)e^{-\ZSI\left ((s-\nu)+(t-\zt)\right )}\right  )Y(\nu,\tii)\ZD\nu
  \right ]\ZD\xii\, \ZD\tii+ \\
\ZLA{eq:SUtrapezio}  &+
\frac{1}{2}
\int _{D(s,t)} 
e^{-\ZSI(s+t)}G(\xii,\tii)
 \ZD\xii\,\ZD\tii\,.
 \end{align} 

 It is sufficient that we prove unique solvability of this equation for at least one value of $\ZSI$.

 Note that if we can prove the existence of a square integrable solution $Y$ then the solution turns out to be continuous (when $f$, $g\in H^1$), thanks to the integrations on the right hand sides. So we prove the existence of a unique square integrable solution of Eq.~(\ref{eq:SUtrapezio}). We need to be precise on the domain of integration. When $T>0$ is fixed, the domain of integration is the set $D(s,t)$ in the figure. Note that the shape of the domain $D(s,t)$ changes according to whether $s\geq t$ or $s<t$ but it is always contained in the large trapezium, which we denote ${\mathcal T}_T$, which is contained in the rectangle   $\mathcal{R}_T=(0,2T)\times(0,T)$.
 
  \begin{center}
\includegraphics[width=10cm]{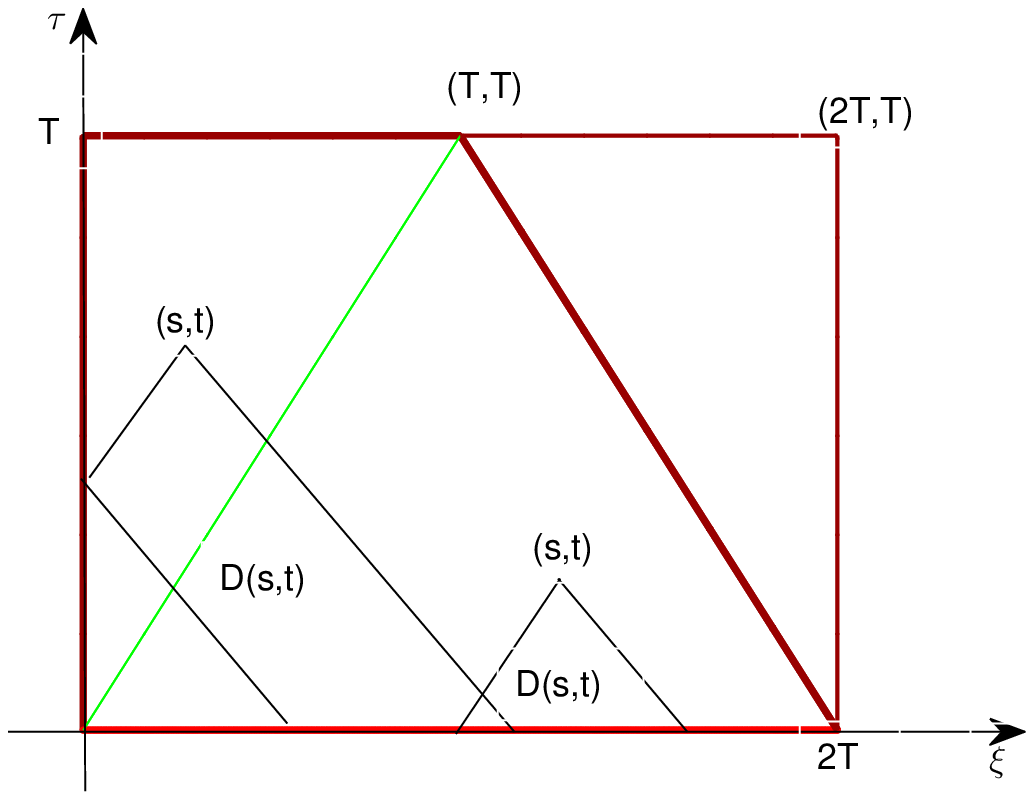} 
   \end{center}
  We consider Eq.~(\ref{eq:intePERblago}) in $L^2(\mathcal{T}_T)$.
We prove that the operator represented by the integral at the first two lines of~(\ref{eq:SUtrapezio}) is a contraction for $\ZSI$ large enough so that the equation has a unique solution.
 This operator is the difference  of the following two operators on $L^2(\mathcal{T}_T)$:
 \begin{align*}
 Y\mapsto L_1(Y)&= \frac{1}{2} \int _{D(s,t)} \left [ 
  \int_0^{\tii}\left ( R_2(\tii-\nu)e^{-\ZSI\left ((s-\xi)+(t-\nu)\right )}\right)Y(\xii,\nu)\ZD\nu \right]\ZD\xii\, \ZD\tii\,,\\
    Y\mapsto L_2(Y)&=\frac{1}{2} \left [\int _{D(s,t)}    \int_0^{\xii}\left (R_2(\xii-\nu)e^{-\ZSI\left ((s-\nu)+(t-\zt)\right )}\right  )Y(\nu,\tii)\ZD\nu
  \right ]\ZD\xii\, \ZD\tii\,.
 \end{align*}
 
 We use boundedness of $R_2(t)$,  $|R_2(t)|<M$, and we prove that for $\ZSI$ large the norms of both these operators are less then $1/2$. We consider the operator $L_1$ (the operator $L_2$ can be handled in a similar way).  In order to simplify the computations, we extend $Y$ with $0$ to $\zzr^2$ so that
 \[
\|Y\|_{L^2 (\mathcal{T}_T)}=\|Y\|_{L^2(K)} \quad \mbox{for every $K\supseteq \mathcal{T}_T$.} 
 \]
Furthermore we note that 
\[
\left \|  \left (L_1Y\right )  \right \|_{L^2({\mathcal{T}_T})}\leq 
\left \|   \left (L_1Y\right )  \right \|_{L^2({\mathcal{R}_T})}  \,.
\]
Note that when $(s,t)\in \mathcal{R}_T\setminus\mathcal{T}_T$ then $D(s,t)\not\subseteq \mathcal{R}_T$ but the contribution to the integral of $Y_{|_{D(s,t)\setminus \mathcal{R}_T}}$  
is zero since $Y$ is zero on this set.
Now we proceed as follows:
\begin{align*}
&\left \|  \left (L_1Y\right )  \right \|^2_{{\mathcal{T}_T}}\leq 
\left \|   \left (L_1Y\right )  \right \|^2_{{\mathcal{R}_T}}=\\
&=\int _{\mathcal{R}_T}\left |
\int _{D(s,t)}\int_0^\zt  \left [ R_2(\zt-\nu)e^{-\ZSI\left ((s-\xi)+(t-\nu)\right )}  \right ] Y(\xi,\nu)\ZD \nu\,\ZD \xi\,\ZD\zt
\right |^2\ZD s\,\ZD t\leq\\
&\leq
 M^2  \int_{\mathcal{R}_T}\left (\left [ \int _{D(s,t)}\int_0^\zt e^{-2\ZSI\left ( (s-\xi)+(t- \nu)\right )}\ZD \nu\,\ZD \xi\,\ZD\zt  \right ]\cdot \right.\\   
&\left.\qquad \cdot\left [
 \int _{D(s,t)}\int_0^{\zt} Y^2(\xi,\nu)\ZD \nu\,\ZD\xi\,\ZD\zt
 \right ]\right)\ZD s\,\ZD t
\end{align*}
It is easy to see that
\begin{align*}
\int_{D(s,t)}\int_0^{\zt} e^{-2\ZSI(s-\xi) }e^{-2\ZSI(t-\nu)}\ZD \nu\,\ZD \xi\,\ZD \zt\leq \frac{1}{2\ZSI}\intt \int _{|s-t+\zt|}^{ s+t-\zt} e^{-2\ZSI(s-\xi)} e^{-2\ZSI (t-\zt)}\ZD\xi\,\ZD \zt=\\
= \frac{1}{4\ZSI^2}\intt e^{-2\ZSI(t-\zt)}\left [
e^{-2\ZSI(\zt-t)}-e^{-2\ZSI\left (s-|s-t+\zt|\right )}
\right ]\ZD\zt \leq \frac{T}{4\ZSI^2}
\end{align*}
(note that we used $s-| s-(t-\zt)|\geq 0$).
Furthermore, using $Y=0$ on $\zzr^2-\mathcal{T}_T$,
\begin{align*}
\int _{D(s,t)}\int_0^\zt Y^2(\xi,\nu)\ZD \nu\,\ZD\xi\,\ZD\zt\leq
\int _{\mathcal{R}_T}\intT Y^2(\xi,\nu)\ZD \nu\,\ZD \xi\,\ZD \zt\leq  T  \|Y\|^2_{L^2(\mathcal{T}_T)}\,.
\end{align*}
We combine these inequalities and we see that
\[
\| L_1 Y  \|^2_{L^2 (\mathcal{T}_T)}\leq \frac{T^4}{2\ZSI^2} 
\|  Y  \|^2_{L^2 (\mathcal{T}_T)}\quad {\rm i.e.}\quad 
\| L_1    \|^2_{\mathcal{L}\left (L^2 (\mathcal{T}_T)\right )}\leq \frac{T^4}{2\ZSI^2} \,.
\]
This inequality shows that $L_1$ has norm strictly less then $1/2$ for $\ZSI$ large enough.

A similar result holds for the operator $L_2$ and this proves the theorem.\zdia
\begin{Remark}
{\rm
We note:
\begin{itemize}
\item we stressed the fact that the differentiations in the previous computations requires smooth boundary term $f$ and $g$. The final integral equation defines a solution for less smooth inputs, as noted in the theorem for the sake of completeness, 
but we are interested in the construction of the quadratic form $\Hfg(T,T)$ and this quadratic form is continuous on $L^2(0,T)\times L^2(0,T)$ for every $T>0$. So, we need not bother with relaxing the regularity assumptions on $f$ and $g$ and we can even work with $f$, $g$ in $\mathcal{D}(0,T)$.

\item the previous arguments in no way depend on the fact that we are considering a   string. Even if $x\in\Omega$, e region with smooth boundary---in this case the integral in~(\ref{eq:defiHfg}) is extended to $\Omega$---the function $\Hfg$ depends on two positive variables and the computations above can be repeated without any change.\zdia
\end{itemize}
}
\end{Remark}

\section{\ZLA{sect:identifQ}Identification of $q(x)$}

 Relaying on the results we obtained up to know, we can show that the   the   method for the identification of $q(x)$ in~\cite{BELIwaves} (see also~\cite{AvdoninBeliIVANOV,AVdonPANDident}) can be adapted to Eq.~(\ref{eq:modello}).
 
We fix a final time  $T\in(0,L]$. Then every $\xi\in L^2(0,T_X)$ (recall the notations: $T_X=T$) is reachable by a control $f$ (see Corollary~\ref{CoroREGOL}). As our target on $(0,T_X)\subseteq(0,L]$ we consider  the solution $\xi(x)$ of 
\begin{equation}\ZLA{eq:statio}
\xi''=q(x)\xi\,,\qquad \xi(0)=0\,,\quad \xi'(0)=1\,.
\end{equation}

For $T$ fixed in  $  (0,L]$, let $f$ be the  \emph{steering control to the target $\xi$,} i.e. the control such that     $w^f(x,T)=\xi(x)$ for $x\in (0,T_X)$. Note that the equality $w^f(x,T)=\xi(x)$  cannot hold for $x\in (T_X,L)$ since on this interval $w(x,T)=0$ while $\xi(x)\neq 0$. So, \emph{the steering control $f$ does depend on $T$, $f=f^T$,} even if $\xi$ does not depend on $T$. 

The function $\xi(x)$ cannot be computed explicitly, since $q(x)$ is unknown but controllability (and $T\leq L$) implies that the steering control exists and it is unique. 
It is an important fact that we can compute this control $f=f^T$ for every $T$ from the response operator (via \Blag equation) even if $q(x)$, hence $\xi(x)$ are unknown. This follows from the following characterization of   the steering control. 
In the derivation of the characterization we use:
\begin{itemize}
\item $\xi(x)$ is smooth, even of class $C^\ZIN$. This justifies the integration by parts below and implies that the steering control   is of class $H^1(0,T)$. So, $y^f(t)\in L^2(0,T)$ (see Theorem~\ref{TeoSEC2RegoLAR}).
 
\item finite velocity of propagation (velocity equal $1$) and  $T\leq L$  
implies that
\[
H^{f,g}(s,t)=\int_0^L w^f(x,t)w^g(x,s)\ZD x=\int_0^{T}w^f(x,t)w^g(x,s)\ZD x\,.
\]
  
\end{itemize}

Let $g\in \mathcal{D}(0,T)$, $T\leq L$. Then we have ($T_X=T$ as usual):

\begin{align}
\nonumber &\Hfg(T,T) = \int_0^L w^f(x,T)w^g(x,T)\ZD x=\int_0^{T_X}w^f(x,T)w^g(x,T)\ZD x=\\
\ZLA{eq:sec4Preli}& 
=\int_0^{T_X} w^f(x,T)\intT w_{,2}^g(x,\zt)\ZD\zt\,\ZD x\,.
\end{align}
 The goal is to give a condition under which $w^f(x,T)=\xi(x)\in C^\ZIN(0,T_X)$. 

We introduce
\[
M(t)=\intt N(r)\ZD r
\]
Let   $f$ be any $H^1$ control   with the additional property    that $w^f(x,T)=w^f(x,T_X)$ is smooth (say of class $H^2$). This is clearly true for the control 
steering to the target $\xi$ defined by~(\ref{eq:statio})). In~(\ref{eq:sec4Preli}),  we replace $w^g_{,2}(x,\zt)$ with its expression~(\ref{eq:modello}) 
and we compute as follows:

 \begin{align*}
&\Hfg(T,T)  =\int_0^{T_X} w^f(x,T)\int_0^{T} \int_0^{\zt} N(\zt-r)\left [
w^g_{,11}(x,r)+q(x)w^g(x,r)
\right ]\ZD r\,\ZD\zt\,\ZD x=\\
&=\int_0^{T_X} q(x)w^f(x,T)\int_0^{T}\left [\int_r^{T}N(\zt-r)\ZD\zt \right ]w^g(x,r)\ZD r\,\ZD x+\\
&+ \int_0^{T_X} w^f(x,T)\int_0^{T}\left [\int_r^{T} N(\zt-r)\ZD \zt\right ]w^g_{,11}(x,r)\ZD r\, \ZD x=\\
&=\int_0^{T }M(T-r)\int_0^{T_X}\left ( q(x)w^f(x,T) +w^f_{,11}(x,T)\right ) w^g(x,r)\ZD x\,\ZD r+\\
&+ \int_0^{T} M(T-r)  \left [w^f(T_X^-,T) w_{,1}^g(  T_X^-,r)- 
   w^f(0,T) w^{g}_{,1}(0,r)-\right.\\
 &  \left.-w^f_{,1}(T_X^-,T)w^{g}(  T_X^-,r)+w^f_{,1}(0,T)w^g(0,r)\right] \ZD r\,.
\end{align*}
 
We are working with $g\in \mathcal{D}(0,T)$ and so (see Corollary~\ref{CoroREGOL})  
$w^g(T_X^{-},r)=0$, $w^g_{,1}(T_X^{-},r)=0$. Hence,  the following equality holds:
 \begin{align*}
&\Hfg(T,T)  =\intT M(T-r)\int_0^{T_X}\left ( q(x)w^f(x,T) +w^f_{,11}(x,T)\right ) w^g(x,r)\ZD x\,\ZD r+\\
&+\intT M(T-r)\left [
-w^f(0,T)y^g(r) +w^f_{,1}(0,T)  g( r)\right] \ZD r\,.
\end{align*}
So we have that $w^f(x,T)=\xi(x)$  solution of~(\ref{eq:statio})  
if and only if the following equality holds:
\[
\langle\Lambda_T^*\Lambda_T f,g\rangle_{L^2(0,T)}
=\Hfg(T,T) =\intT M(T-r)g(r)\ZD r\,.
\]
Using controllability, i.e. invertibility of   $\Lambda_T^*\Lambda_T$,
we compute $f(t)$, the steering control to the \emph{unknown} target $\xi(x)$, $x\in (0,T_X)$, when $\xi(x)$ is given  given by~(\ref{eq:statio}):
\[
f(t)=f^T(t)=\left [ \left (\Lambda_T^*\Lambda_T\right )^{-1}M(T-\cdot)\right ] (t)\,.
\]
We repeat that $\Lambda_T^*\Lambda_T$ is known,  computed from \Blag equation.

Of course, $w^f(x,t)$, $t\in (0,T)$, cannot be computed since $q(x)$ is still uknown. Nevertheless, $\xi(T_X)=\xi(T)$ can be computed since, from Item~\ref{TeoSEC2RegoLARitem1BIS} of Theorem~\ref{TeoSEC2RegoLAR},
we have
\[
\xi(T_X)=w^f(T_X^-,T)=f(0)=f^T(0)\,.
\]
We repeat  (in principle!) this computation for every $T\in (0,L)$ and we compute the function $\xi(T)$ for every $T\in (0,L)$ from which we get
\[
q(T)=\frac{\xi''(T)}{\xi(T)}
\]
(extension by continuity at the finitely many points in which $\xi =0$ and to $T=0$).

This is the extension to system~(\ref{eq:modello}) of the algorithm presented in~\cite{AvdoninBeliIVANOV,BELIwaves}.

\section{Appendix: the proofs of the results of Sect.~\ref{sectOBERCmodel}}

In this section we prove theorems~\ref{Teo:RegoSOLUZ} and~\ref{TeoSEC2RegoLAR}.

\subparagraph*{The proof of Theorem~\ref{Teo:RegoSOLUZ}}
We can prove the regularity for the function $W(x,t)$ which is obtained after the MacCamy trick.  Thanks to the fact that $T\leq L$, the solution $W(x,t)$ can be represented by using  Formula~(\ref{FormDALA}) with
\begin{equation}\ZLA{eq:LaFinAPPeN}
F(x,t)=\left ( q(x)+\zaa\right ) W(x,t)+\intt K(t-s) W(x,s)\ZD s\,,
\end{equation}
i.e.
\begin{equation}\ZLA{eq:LaFinAPPeNWu}
  \left\{\begin{array}{l}W(x,t) =e^{-\zg(t-x)}f(t-x)+\frac{1}{2}u(x,t)\,,\\[3mm]
u(x,t) =\intt \int_{|x-t+\zt|}^{x+t-\zt}  F(\xi,\zt)\ZD\xi\,\ZD \zt=\\[2mm]
 =
\intt \int_{|x-t+\zt|}^{x+t-\zt}\left [\left ( q(x)+\zaa \right )W(\xi,\zt)+\int_0^\zt K(\zt-s) W(\xi,s)\ZD s\right ]\ZD\xi\,\ZD\zt\,. 
\end{array}\right.
\end{equation}
By assumption, $f\in\mathcal{D}(0,T)$ and so $u$ and $W$ have the same regularity and we know $W\in C([0,T]L^2(0,L))$ so that $u$ is (at least) of class $H^1$.
Note that

\begin{align*}
u(x,t)&= \int_0^{t-x} \int _{t-x-\zt}^{x+t-\zt} F(\xi,\zt)\ZD\xi\,\ZD\zt+\int _{t-x}^t \int _{x-t+\zt}^{x+t-\zt} F(\xi,\zt)\ZD\xi\,\ZD\zt\,,\quad 0<x<t\\
u(x,t)&=
\int _0^t \int _{x-t+\zt}^{x+t-\zt} F(\xi,\zt)\ZD\xi\,\ZD\zt\,,\quad t<x<L
\,.
 \end{align*}
 It follows that $u\in C([0,L]\times[0,T])$ and so we have also $F\in C([0,L]\times[0,T])$.
  Computing the derivatives we see that 
\[\begin{array}{lll}
&  \mbox{
\hskip 3cm when  $ 0<x<t$  
} &\\[1.5mm]
u_t(x,t)&=\intt F(x+t-\zt,\zt)\ZD\zt+\int _{t-x}^t F(x-t+\zt,\zt)\ZD\zt-\int_0^{t-x} F(t-x-\zt,\zt)\ZD\zt=&\\[1mm]
&=\int_x^{x+t} F(\xi,x+t-\xi)\ZD\xi+\int_0^x F(\xi,\xi-x+t)\ZD\xi-\int_0^{t-x}F(\xi,t-x-\xi)\ZD\xi\,,&\\[1mm]
u_x(x,t)&=\int_x^{x+t} F(\xi,x+t-\xi)\ZD\xi-\int_0^x F(\xi,\xi-x+t)\ZD\xi+\int_0^{t-x} F(\xi,t-x-\xi)\ZD\xi\,.&
\end{array}\]
Analogously 
\[
\begin{array}{lll}
&\mbox{\hskip 3cm when $t<x<L$ }&\\[1.5mm]
u_t(x,t)&=\int_x^{x+t}F(\xi,x+t-\xi)\ZD\xi+\int _{x-t}^x F(\xi,\xi-x+t)\ZD\xi\,,&\\[1mm]
 u_x(x,t)&=\int_x^{x+t}F(\xi,x+t-\xi)\ZD\xi-\int _{x-t}^x F(\xi,\xi-x+t)\ZD\xi &
\end{array}
\]
(we might note, from the expression~(\ref{eq:LaFinAPPeN}) and finite velocity of propagation, that $F(\xi,\zt)=0$ when $\xi>\zt$, but we don't need the use of this piece of information here).

We note from these expressions that $u_x(x,t)$ and $u_t(x,t)$ are both continuous on $[0,L]\times [0,T]$ (the diagonal $x=t$ included). Hence, $W(x,t)\in C^1([0,L]\times [0,T])$ and, from~(\ref{eq:LaFinAPPeN}), $F(x,t)\in C^1([0,L]\times [0,T])$.

Now we compute the second derivatives, both for $0<x<t$ and for $t<x<L$ and we see that they are continuous on these sets, and we see that $W$ satisfies the equality~(\ref{eq:modello}) on these sets.

We pass to the limit for $x\to t^+$ and $x\to t^-$ and we use $F(0,0)=0$. We see that the partial derivatives are continuous and  that equality~(\ref{eq:modello}) holds on $[0,L]\times[0,T]$.
Note that the condition $F(0,0)=0$ follows from~(\ref{eq:LaFinAPPeN}) since $W(x,t)$ is continuous on $[0,L]\times[0,T]$ and $W(x,0)=0$ (analogously, $W(0,t)=0$ for $t$ ``small'', since $f\in\mathcal{D}(0,T)$).\zdia

\subparagraph*{The proof of Theorem~\ref{TeoSEC2RegoLAR}}
We apply MacCamy trick and we prove the result using Eq.~(\ref{eq:DOPOmaCCamy}) of $W(x,t)$.  So, $W(x,t)$ is given by~(\ref{eq:LaFinAPPeN})-(\ref{eq:LaFinAPPeNWu}) since the standing assumption in this theorem is $T\leq L$.

\emph{Proof of  the statement in item~\ref{TeoSEC2RegoLARitem1}.}
  We know from~\cite{PandCINA,PandLIBRO} that $W\in C([0,T];L^2(0,L))$ (zero on $(T,L)$). Hence, the function $u(x,t)$ is a continuous function of $(x,t)$: \emph{the discontinuity of $W$, hence of $w$, are the same as those of $e^{-\zg(t-x)}f(t-x)$.}

  Let $w^f(x,T)=\xi(x)$. Then,
  \begin{equation}\ZLA{eq:ExpreSTERinTargREG}
 \xi(x)=e^{ \zg x  }f(T-x)+\frac{1}{2} u(x,T)\,.
  \end{equation}
  If $\xi(x)$ is continuous, then $f(t)$ is continuous too. Furthermore, the kernel of the transformation $W\mapsto u$ is zero (even in the case $q(x)+ \zaa=0$) so that $f(t)$ is uniquely identified.
  
 This is the statement in item~\ref{TeoSEC2RegoLARitem1}. 
 
 \emph{We prove the statement in item~\ref{TeoSEC2RegoLARitem1BIS}.}
We use finite velocity of propagation, so that $w(x,t)=0$ if $x>t$ and $\lim _{x\to t^+ } w(x,t)=0$. Note that we have also $F(x,t)=0$ if $x>t$ and the definition of $u(x,t)$ shows that
\[
\lim _{x\to T^-} u(x,T)=0\,.
\]  
  In particular, using the fact that $w$ and $W$ share the same regularity properties, we see that the limit $\lim _{x\to T^-}w(x,T) $  exists if and only if $\lim _{t\to 0^+} f(t)$ exists and:
  \[\ZLA{eq:delLIMITaZERO}
 \lim _{x\to T ^-} w(x,T )=\lim _{t\to 0^+} f(t)\,,\qquad \lim _{x\to t ^+} w(x,T )=0\,, \quad \mbox{when $T\leq L$}\,.  
  \]
 This implies the property in the statement~\ref{TeoSEC2RegoLARitem1BIS} since $f(x)$ is continuous when $\xi(x)$ is continuous. 

\emph{Now we prove the property stated in item~\ref{TeoSEC2RegoLARitem2}.  }
The derivative $u_x(x,t)$ has been computed in the proof of Theorem~\ref{Teo:RegoSOLUZ}. Continuity of the shift in $L^2$ proves that
the following limit exists in $L^2(0,T)$:
\[
u_{x}(0,t)=\lim _{x\to 0^+}u_{x}(x,t) = \int_0^tF(\xi,t-\xi)\ZD\xi\,.
\]
 This is a continuous function of $t$ for every square integrable boundary control $f$ and so
 $w_x(0,t)$ exists in $L^2(0,T)$ if and only if $f$ is of class $H^1$.

Finally, let $w^f(x,T)=\xi(x)\in H^1(0,T_X)$, $T_X=T\leq L$. Then $\xi(x)$ is continuous and $f(t)$ is continuous too, from the statement in item~\ref{TeoSEC2RegoLARitem1}. Hence, 
$u_x(x,T)\in L^2(0,T_X)$ (see the proof of Theorem~\ref{Teo:RegoSOLUZ}). Formula~(\ref{eq:ExpreSTERinTargREG}) shows that 
$f\in H^1(0,T)$. \emph{This is the proof of  item~\ref{TeoSEC2RegoLARitem4}.}\zdia

\end{document}